\documentclass[12pt]{article}
\usepackage{amsmath}\usepackage{psfig}
\usepackage{amssymb}
\begin{document}
\title{The geometry of bifurcation surfaces in parameter space. I. 
A walk through the pitchfork}
\author{Rowena Ball \\
Department of Theoretical Physics\\ 
Research School of Phys. Sci. \& Eng.\\ The Australian National University\\
Canberra 0200 Australia\thanks{e-mail: rxb105@rsphysse.anu.edu.au}}
\date{}
\maketitle
\begin{abstract}
The classical pitchfork of singularity theory is a twice-degenerate 
bifurcation that typically occurs in dynamical system models exhibiting 
$\mathbb{Z}_2$ symmetry. Non-classical pitchfork singularities also occur
in many non-symmetric systems, where the total bifurcation environment is
usually more complex. In this paper three-dimensional manifolds of 
critical points, or limit-point shells, are introduced by examining 
several bifurcation problems that contain a pitchfork as an organizing
centre. Comparison of these surfaces shows that notionally equivalent 
problems can have significant positional differences in their bifurcation
behaviour. As a consequence, the 
parameter range of jump, hysteresis, or 
phase transition phenomena in dynamical models (and the physical systems 
they purport to represent) is determined by other singularities that
shape the limit-point shell.
\end{abstract}
\begin{center}{\sl Keywords: pitchfork; singularity theory; bifurcation theory; 
limit-point shell\\
PACS classification numbers: 47.20.Ky, 05.45.-a}
\end{center}
\renewcommand{\baselinestretch}{1.1}\normalsize
\section{\label{intro}Introduction}
Discontinuous behaviour in dissipative dynamical systems is often 
ascribed to the propinquity of a pitchfork, 
a codimension 2 singularity which requires two auxiliary parameters for a 
universal unfolding. In a 3-dimensional space labelled by the system 
parameters a universal unfolding forms a critical surface, or limit-point shell.
{\it A priori} knowledge of the topology 
of critical surfaces, and how they change with variation of other 
parameters (i.e., animation) is a powerful aid in the design 
and optimization of dynamical systems and the control of jump, 
hysteresis, or oscillatory phenomena.  In this work cartoons of computed 
bifurcation surfaces 
are used to visualize the pitchfork and its surroundings in parameter space. 
It is shown how these surfaces are shaped by the degree of symmetry and by
the lower codimension antecedents of the pitchfork as the organizing 
centre, or highest order singularity,
in particular bifurcation problems.\footnote{``Bifurcation''
simply means ``fork'', so it seems unnecessary to refer to a pitchfork fork.}

A mathematical description of the pitchfork in singularity theory terms 
will be given in \S \ref{s2}, for now the name itself will suffice as
a working definition. The classic trident representation is simply a graph
depicting how the number of something (such as solutions of an 
equation, or states of an ideal ensemble) changes from one to three as 
a control variable is made to cross a critical value. Because of both its
heuristic value in singularity theory and its importance in applied 
bifurcation problems, the pitchfork has been described in many texts 
dealing with dynamics and bifurcations. A recent exposition with this flavour 
can be found in \cite{Kuznetsov:1998}. 

Pitchforks have been most often reported as bifurcations of the steady
states of idealized systems possessing $\mathbb{Z}_2$ or reflectional 
symmetry. Simple examples in this category include studies of the driven 
pendulum 
\cite{Broer:1999}, \cite{Schmitt:1998}, \cite{Johnson:1998}, and 
variations on the buckling beam problem \cite{Yabuno:1998}, 
\cite{Golubitsky:1979}. In these systems, which are also called 
imperfect bifurcation problems, symmetry-breaking perturbations dissolve
the pitchfork leaving behind persistent limit points where discontinuous
action may occur. In more complicated problems the symmetric pitchfork occurs
as a bifurcation of periodic solutions \cite{Jensen:1999}, or may be 
embedded in a higher-order degeneracy \cite{Algaba:1999}.  

Although the presence of a pitchfork in a mathematical model can
sometimes be interpreted as {\it prima facie} 
evidence or diagnostic of $\mathbb{Z}_2$ invariance, there are also
many bifurcation problems containing non-symmetric pitchforks.  In 
some problems of this type a fully unfolded or perturbed pitchfork is intrinsic
to a minimal description of the associated physical system. Thermokinetic and 
isothermal chemical systems are the most widely studied in this category; a 
good introduction to these is found in \cite{Gray:1990}.
A partially unfolded pitchfork occurs in other non-symmetric problems, with a 
single static perturbation required to complete a realistic description of the
system. Two recent examples are a model for L(low)--H(high) confinement state
transitions in plasmas \cite{Sugama:1995,Ball:1999c}
and experimental and numerical 
studies of an electronic Van der Pol oscillator \cite{Juel:1997}.

With recent advances in inertial manifold
theory \cite{Stewart:1999}, it is likely that the pitchfork in general will
also appear in the dynamics of infinite-dimensional systems that exhibit
low-dimensional behaviour on a long time-scale. 

In this paper several bifurcation problems will be 
discussed and compared to exemplify the qualities of pitchfork manifolds in the
different categories: (1) the prototypic universal unfolding of the 
pitchfork, (2) a well-known thermokinetic system, (3) the simplest universal
unfolding of a non-symmetric pitchfork, and (4) the L--H transition problem
mentioned above.

Bifurcation surfaces can also illuminate rather dramatically the truism that 
what you see depends on where you view the object from. Originally the
pitchfork was described from an orthogonal point of view as the generic
cusp, and in \S \ref{s2} we work around to the prototypic pitchfork 
beginning from the more familiar cusp manifold. Moving from the non-generic
surface of equilibrium points or steady states, up a dimension to the invariant
surface of fold points, the limit-point shell $L_p$ of the prototypic pitchfork
is presented as a fundamental, generic object. In \S \ref{s3} the 
limit-point shell of the CSTR problem, $L_c$ (first described in 
\cite{Ball:1999a} and \cite{Ball:1999b}), 
and that of the simplest non-symmetric pitchfork, $L_{\scriptscriptstyle TI}$, 
are compared.
Some curious features of the extraordinary limit-point shell of  
the L--H problem, $L_{\scriptscriptstyle LH}$ are also described in 
\S \ref{s3}. 
\S \ref{s4} concludes with a
brief discussion of what these result imply for the design and control of 
experimental dynamical systems.
\section{\label{s2}From cusp manifold to limit-point shell}
\begin{figure}\hspace*{2cm}
\psfig{file=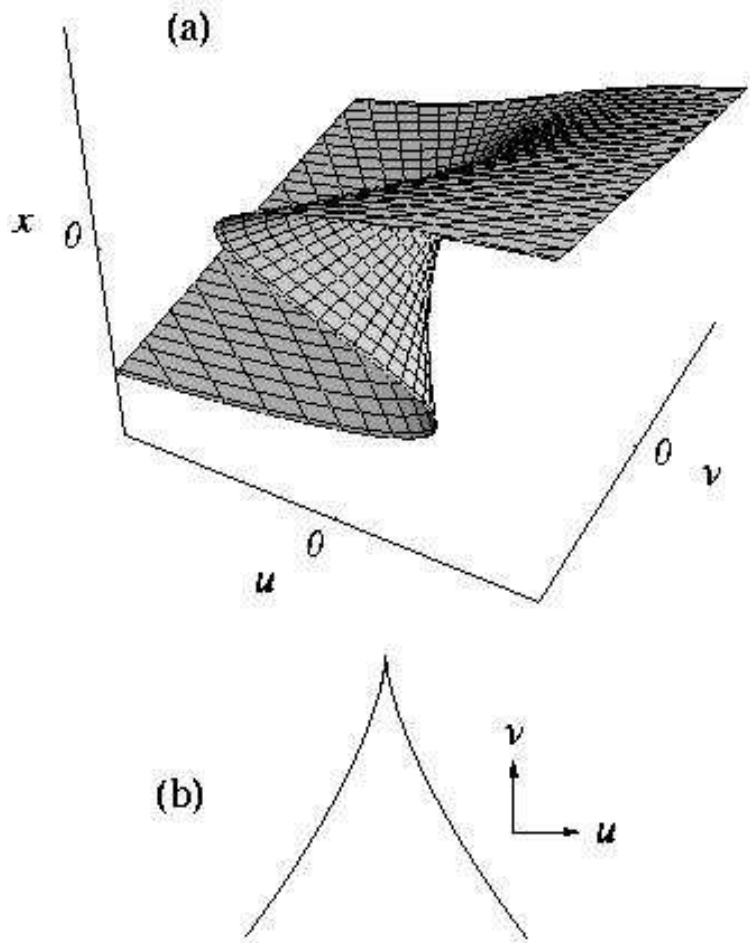,width=10cm}
\caption{\label{fig1}The cusp catastrophe.}
\end{figure}

\begin{figure}\hspace*{2cm}
\psfig{file=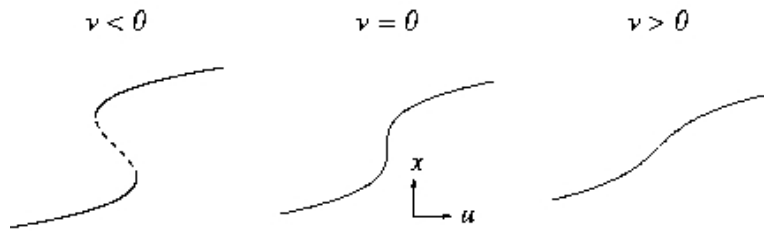,width=10cm}
\caption{\label{fig2}Paths across the cusp manifold.}
\end{figure}

In 1955 Hassler Whitney published the first exposition of singularity theory
\cite{Whitney:1955}, in which he derived conditions for a regular point
$p$ of a smooth mapping $f$ from $\mathbb R^2$ into $\mathbb R^2$ to be a cusp 
point.
In coordinates $(u,v)$, $(x,y)$ these are
\begin{equation}\label{cc}
\begin{aligned}
\text{a.}&\hspace{1cm}u_x=u_y=v_x=0, \, v_y=1\quad\text{(singular condition)},\\
\text{b.}&\hspace{1cm}u_{xx}=0, \, u_{xy}\neq 0, \, 
u_{xxx}-3u_{xy}v_{xx}\neq 0,\quad
\text{(cusp condition)}
\end{aligned}
\end{equation} 
at $p$. Whitney proved that any mapping containing a point at which
the conditions (\ref{cc}) are satisfied can be transformed by coordinate
changes into the following normal form for the cusp:
\begin{equation}\label{nfc}
\begin{aligned}
u & = xy - x^3\\
v & = y
\end{aligned}
\end{equation}
Subsequently Ren\'{e} Thom listed the cusp, which he called the Riemann-Hugoniot 
catastrophe, as the second of the famous (I mean, of course, infamous) seven
elementary catastrophes \cite{Thom:1972}. In catastrophe theory the normal form 
(\ref{nfc}) becomes the universal unfolding $G(x,u,v)$ of the 
germ $g(x)=x^3$:
\begin{equation}\label{uuc}
G(x,u,v)=x^3-vx+u,
\end{equation}
where $G$ is the gradient of a corresponding governing potential $V$. 
The familiar cusp surface, $G(x,u,v)=0$, 
 and the projection of the folds on the $u,v$ plane are
shown in figure \ref{fig1}. 
The surface may be viewed as the lateral unfurling of a path in the 
$x,u$ plane
into $v$. Three qualitatively different paths --- slices of the cusp 
surface at constant $v$ or bifurcation diagrams ---
are sketched in figure \ref{fig2}.

\begin{figure}
\psfig{file=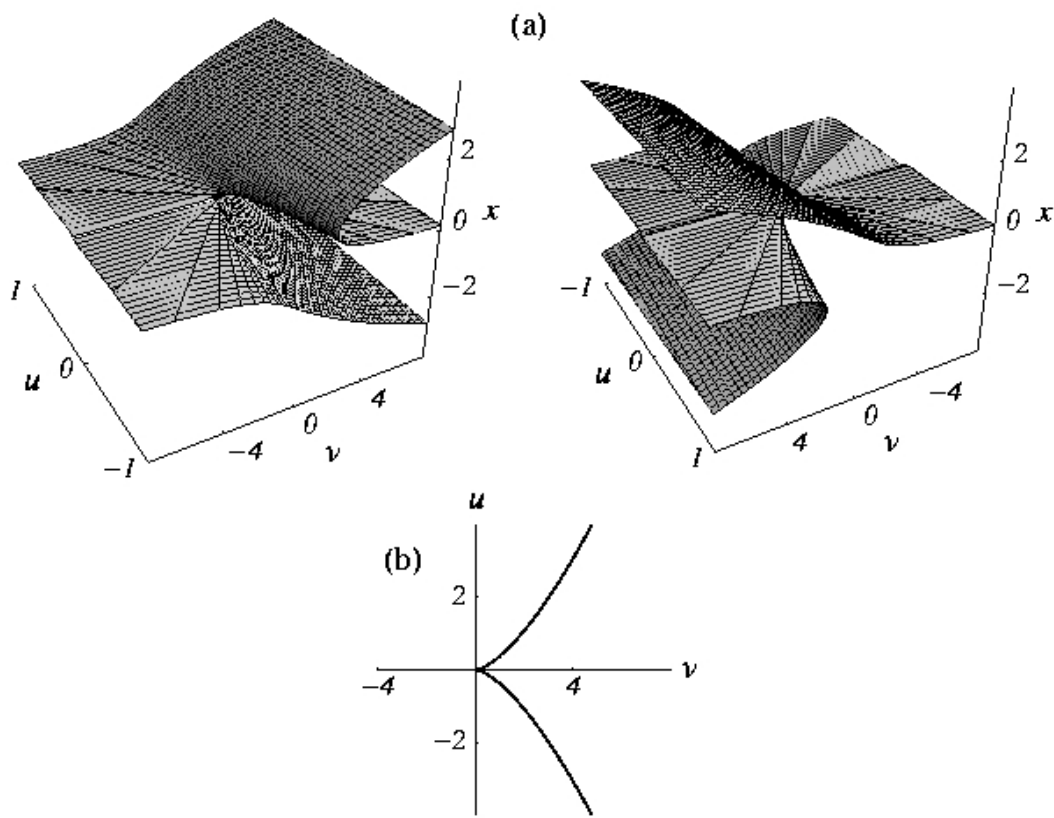,width=14cm}
\caption{\label{fig3}An orthogonal path through the 
cusp unfolding opens into a manifold around the pitchfork.}
\end{figure}

\begin{figure}\hspace*{1.5cm}
\psfig{file=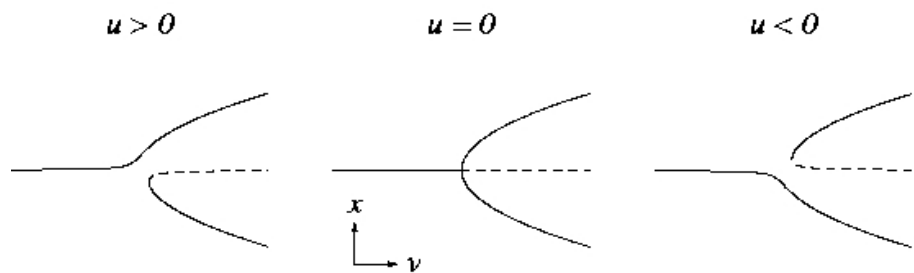,width=12cm}
\caption{\label{fig4}Paths across the pitchfork manifold
in figure \ref{fig3}(a).}
\end{figure}

Although the cusp is generic, in the sense that all other singularities may
be perturbed to either a fold or a cusp \cite{Arnold:1992}, the surface
in figure \ref{fig1} is not
a unique manifold of the cusp catastrophe. Since all paths through the 
cusp unfolding (\ref{uuc}) are equally valid, 
we may choose a path in the $x,v$
plane. Any such path unfurls laterally into $u$ to form another surface, 
shown from two points of view 
in figure \ref{fig3} together with 
the projection of the folds on the $u,v$ plane. 
Three qualitatively different bifurcation diagrams are shown in 
figure \ref{fig4}, whence we finally arrive at the classical  
pitchfork.

Returning to Whitney's original defining conditions for
 the cusp, equations (\ref{cc}), 
 we see now that they also define the pitchfork. 
 
\begin{figure}\hspace*{3cm}
\psfig{file=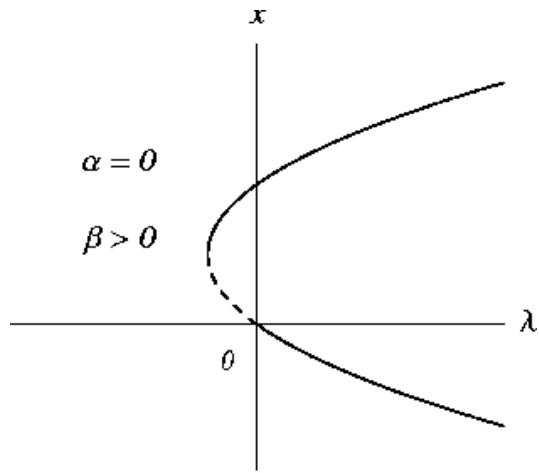,width=7cm}
\caption{\label{fig5}The bifurcation diagram of $P_p$
for $\alpha=0$, $\beta=2$.}
\end{figure}

 The surface in figure \ref{fig3} is non-generic, because it
 does not represent {\em all} qualitative 
 information about the pitchfork. Golubitsky and Schaeffer 
 \cite{Golubitsky:1985} proved 
 that a universal unfolding of the pitchfork must include a fourth variable. 
 Their classification of singularities by codimension and derivation of 
 universal unfoldings has the contextual setting of bifurcations
 of steady states in autonomous dynamical systems dependent on 
 parameters:
 \begin{equation}\notag
 \frac{dx}{dt}=G\left(x,\lambda,\alpha_i\right)=0.
 \end{equation} 
Here $x$ is the dynamical state variable, $\lambda$ is the principal bifurcation
parameter and the $\alpha_i$ are auxiliary or unfolding parameters.
Assuming henceforth this context and notation, the pitchfork conditions 
$P$ are given in table \ref{table1}, column 2.
A bifurcation problem which satisfies $P$ is said to be locally equivalent to
the normal form
 \begin{equation}\label{nfp}
 g\left(x,\lambda\right)=\pm x^3\pm \lambda x.
 \end{equation}  
The prototypic universal unfolding of the pitchfork, $P_p$, is given as
\begin{equation}\tag{$P_p$}
 G\left(x,\lambda,\alpha,\beta\right)=x^3-\lambda x+\alpha+\beta x^2.
 \end{equation} 
 A na\"{\i}ve interpretation of the concept of a universal unfolding is more
 useful here than a rigorous definition or derivation. 
 We apply the conditions $P$ to the unfolding
 \begin{equation}\label{part}
 G\left(x,\lambda,\alpha\right)=x^3-\lambda x +\alpha.
 \end{equation}  
 (This is a partial unfolding of (\ref{nfp}); compare equation (\ref{uuc}).)
 The result comprises an overdetermined system of four equations in 
 three variables. 
 If we perturb (\ref{part}) at $\alpha=0$ 
 by adding the term $\beta x^2$, the bifurcation diagram acquires qualitatively
  different characteristics, as indicated in figure \ref{fig5}. 
However, more perturbations, or different perturbations, do not introduce 
any more qualitative differences to the bifurcation diagram.  
(Note: it can be guessed that universal unfoldings are not necessarily unique.)

\begin{figure}\hspace*{1cm}
\psfig{file=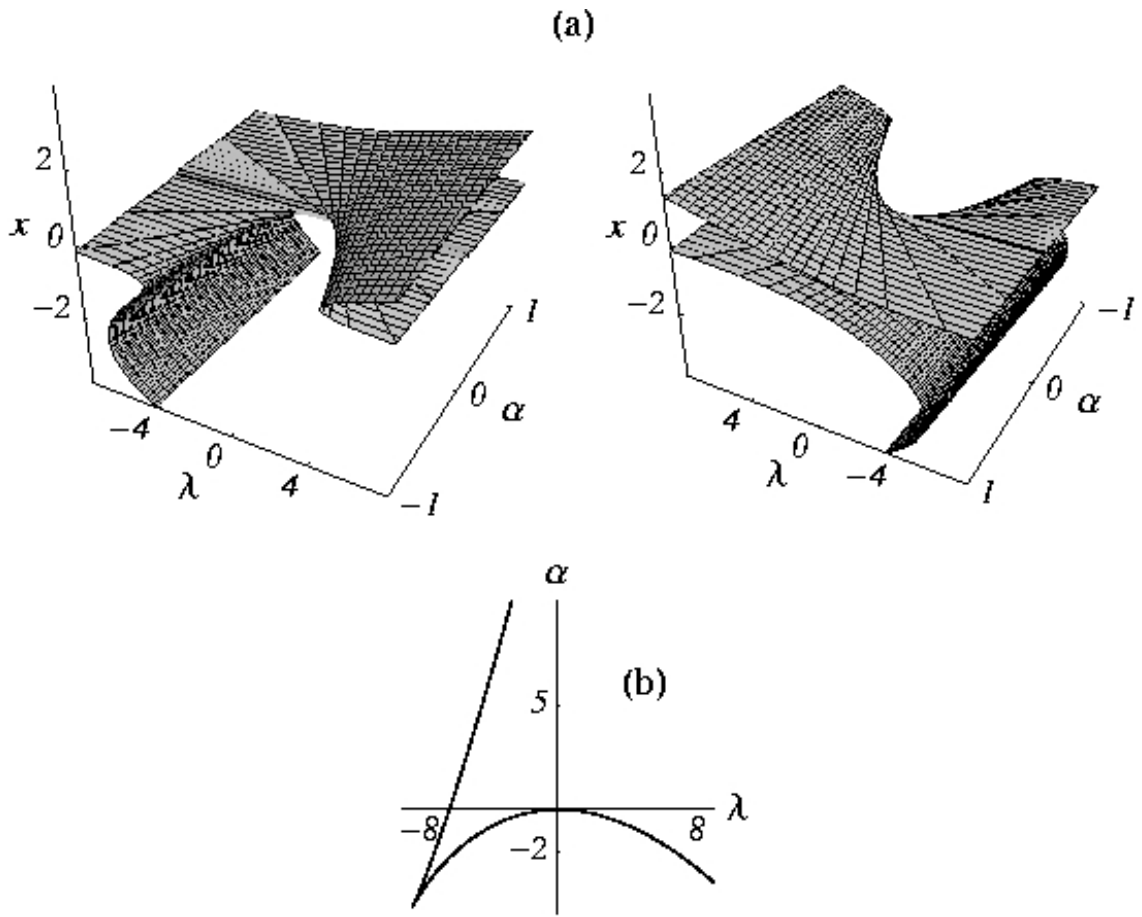,width=15cm}
\caption{\label{fig6}A bifurcation manifold of $P_p$, $\beta=5$.}
\end{figure}

\begin{figure}
\psfig{file=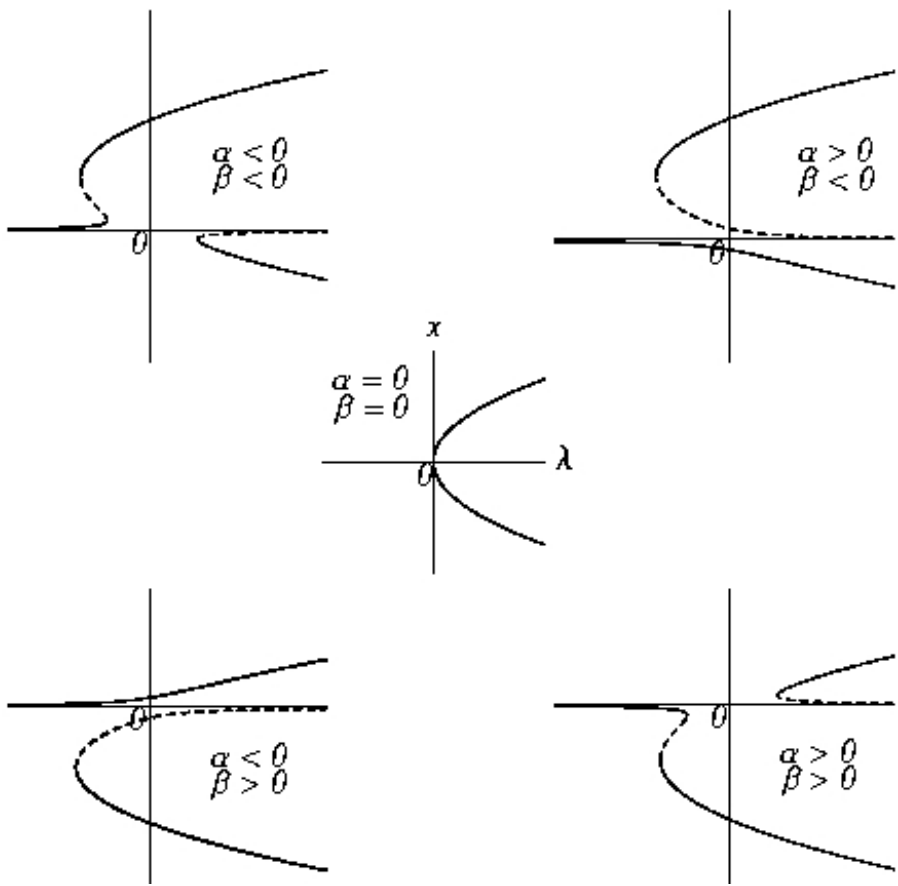,width=12cm}
\caption{\label{fig7}The qualitatively distinct bifurcation diagrams of $P_p$.}
\end{figure}

For any value of $\beta\neq 0$ the bifurcation 
surface of $P_p$ formed by unfolding a path in
the $x,\lambda$ plane into $\alpha$ does not, therefore,
 include the pitchfork. Two views
of this surface for $\beta=5$ are shown in figure \ref{fig6}
along with the projection of the folds on the $\lambda,\alpha$ plane.
The famous five bifurcation diagrams for the universal unfolding 
$P_p$ are given in figure~\ref{fig7}. 

\begin{figure}
\psfig{file=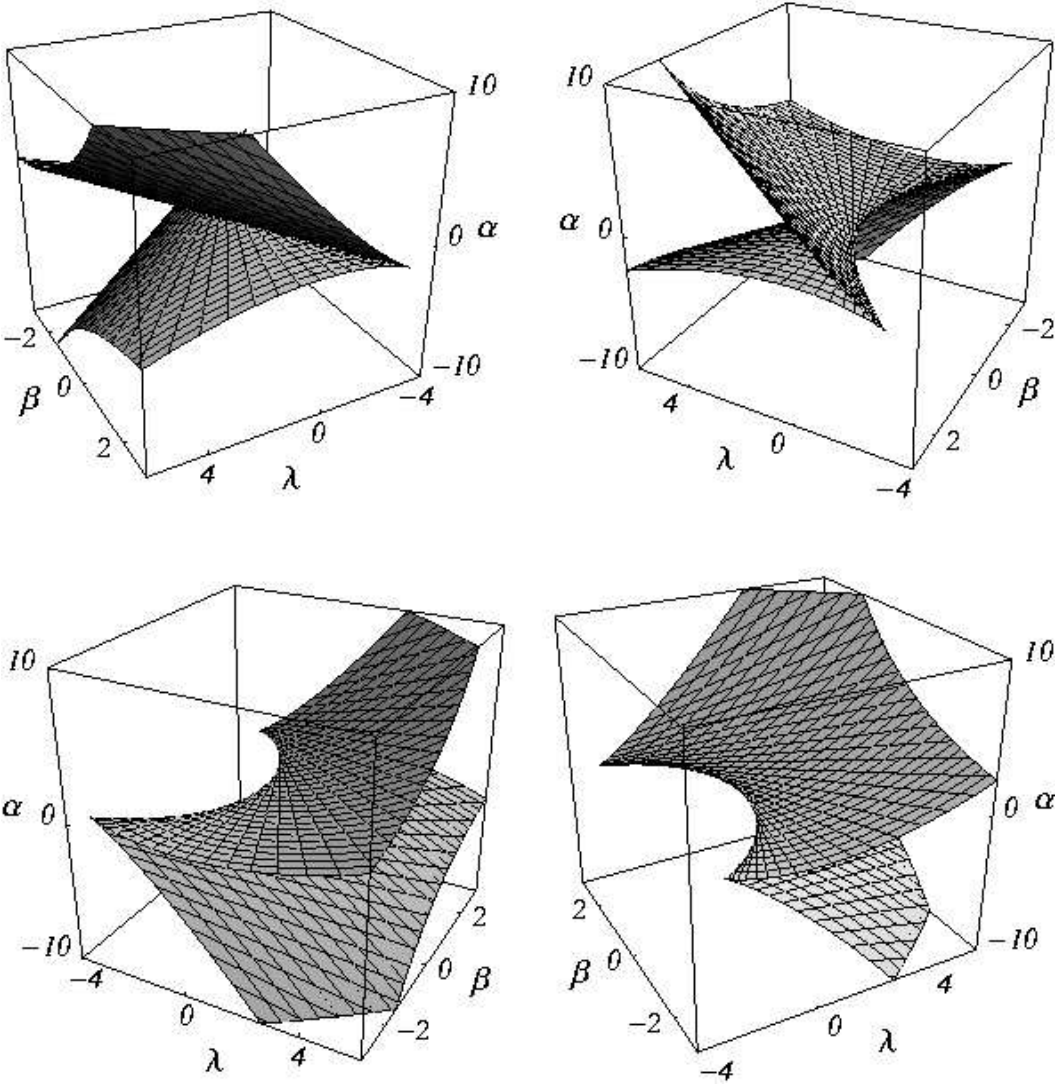,width=14cm}
\caption{\label{fig8}Four views of $L_p$.}
\end{figure}

There are two singularities on the fold lines of figure \ref{fig6}(b), 
a hysteresis point $H$
and a transcritical point $T$, definitions of which are given in 
in table \ref{table1} columns 2 and 3.
It is evident from these definitions and from the surface
for the partial unfolding in figure \ref{fig3} that the pitchfork is the
limiting degeneracy of $H$ and $T$ for $\beta=0$. A unique, continuous 
manifold $L_p$ around $P_p$ is obtained by unfurling the {\em fold} lines in
the $\lambda,\alpha$ plane into $\beta$. This forms a 
limit-point shell, 
 the surface of fold or limit points of an unfolding. 
$L_p$ is shown from four vantage points in figure \ref{fig8}. It is a
 self-contained and unique manifold around $P_p$; it cannot be 
translated, rotated, deformed or punctured. 
Let us inspect it closely. 

A slice of the surface at constant $\beta\neq 0 $ is simply the
fold lines in figure \ref{fig6}(b). There is a distinct seam of hysteresis 
points that runs from positive to negative $\beta$ and negative to 
positive $\alpha$ through the pitchfork 
at $(0,0,0)$. Less obvious is the line of trancritical points, which lies
along $\lambda =0$ in the  $\lambda,\beta$ plane.
 The shell is symmetric about two 
reflections: reflect it in a vertical mirror plane through $\beta=0$, then 
reflect it again in a horizontal mirror plane through $\alpha=0$.
\section{\label{s3}The limit-point shell in applications}
In dynamical models the limit-point shell outlines
 the boundary and extent of steady-state multiplicity over
the parameter space.
Therefore, the shape of a limit point shell is an important consideration
in predictive work, because it is a guide in the selection of design criteria
and operating conditions for an experimental system.
However, $L_p$ is a poor metaphor for many models which contain a
 pitchfork. A specific limit-point shell is shaped by the 
presence and location of other bifurcations and by the symmetry of the
problem. In this section these ideas are illustrated by 
comparing the limit point shells of three bifurcation problems which have
pitchforks as organizing centres. 
\subsection{\label{cstrproblem}The CSTR problem}
Interest has been maintained in this thermokinetic system since the 1950s, 
because it posesses a combination of unusual dynamical properties and 
real-world experimental accessibility. The first singularity theory
study of the CSTR problem can be
found in \cite{Golubitsky:1980} and a recent alternative treatment 
is given in \cite{Ball:1999a}.

\begin{figure}\hspace*{3cm}
\psfig{file=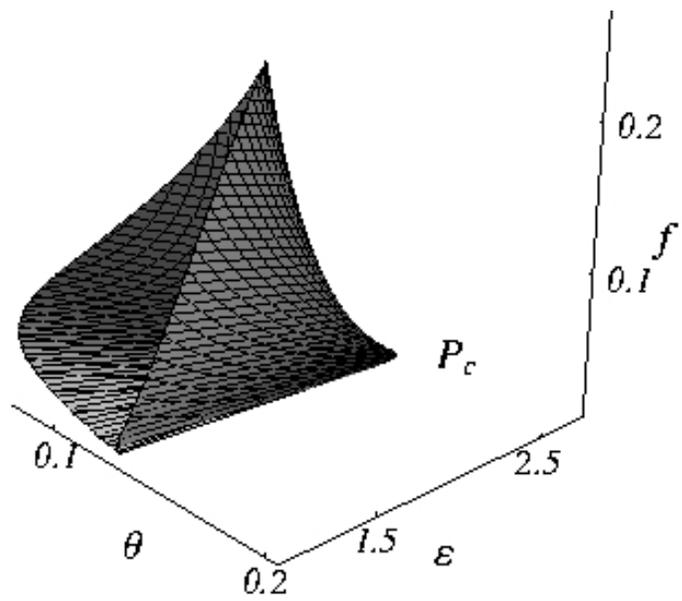,width=9cm}
\caption{\label{fig9}$L_c$ for $\ell=0.05$.}
\end{figure}

The simplest CSTR model describes an exothermal chemical reaction occurring 
in a well-stirred bounded medium. 
As a bifurcation problem it may be written as follows:
\begin{equation}\tag{$P_c$}
G(u,f,\theta,\varepsilon,\ell)=
\frac{f e^{-1/u}}{e^{-1/u}+f}+\left(\varepsilon f+\ell\right)
\left(\theta-u\right).
\end{equation}
The state variable is the temperature $u$
which depends on a number of parameters: $f$, an input rate; $\theta$, the
coupled temperature of the thermostat and the input; $\ell$, the thermal
dissipation rate; and $\varepsilon$, the intrinsic properties of the
medium. For $u\neq 1/2$ it can be shown that an 
organizing centre is the pitchfork $P_c$,  although the system has
no symmetries. The limit point shell $L_c$ at a fixed value of $\ell$ is shown
in figure \ref{fig9}. 
({\it Note:} A detailed discussion of this structure may be found in 
\cite{Ball:1999a}. The dynamical CSTR problem also contains Hopf bifurcations,
the manifold of which will be discussed elsewhere.) 

Contrasted with $L_p$, the prototypic limit point shell
of figure \ref{fig8}, the overall qualities of $L_c$ appear to be 
{\em asymmetry} and {\em convexity}. 
\subsubsection*{\label{lob} The role of co-existing bifurcations}
Table \ref{table1} shows that embedded in the pitchfork are the codimension 1 
bifurcations $H$ and $T$ {\em or} $I$, or $H$ and $T$ {\em and} $I$, since
$T$ and $I$ are distinguished only by the sign of the non-degeneracy condition
$\det d^2G\neq 0$.

For the prototypic unfolding, $P_p$,
 we find that $G=G_x=G_\lambda=0$ and 
$\mathrm{det}\,(d^2G) = -1$, $G_{xx}=-6$ at 
$x=\lambda=\alpha=0$, thus $T$ but not $I$ is embedded in $P_p$ . 
For the CSTR problem it is shown in Appendix A
that both $T$ {\em and} $I$ are embedded in $P_c$. 

The last column in table \ref{table1} gives conditions for another
codimension 2 bifurcation, the asymmetric cusp $A$. In Appendix A it is 
shown that this singularity is also present in the CSTR problem. 
$A$ is too inconspicuous
to be pinpointed visually on 
the portion of $L_c$ shown in figure \ref{fig9}, but it is
easily approximated numerically as $(u,f,\theta,\varepsilon)\approx
(0.216,0.016,0.141,1.91)$ at $\ell=0.05$. Nevertheless, the form of $L_c$ is
strongly sculpted by the presence of $A$.
\subsubsection*{Is the normal form for the CSTR problem an adequate proxy?}
\begin{figure}\hspace*{1.5cm}
\psfig{file=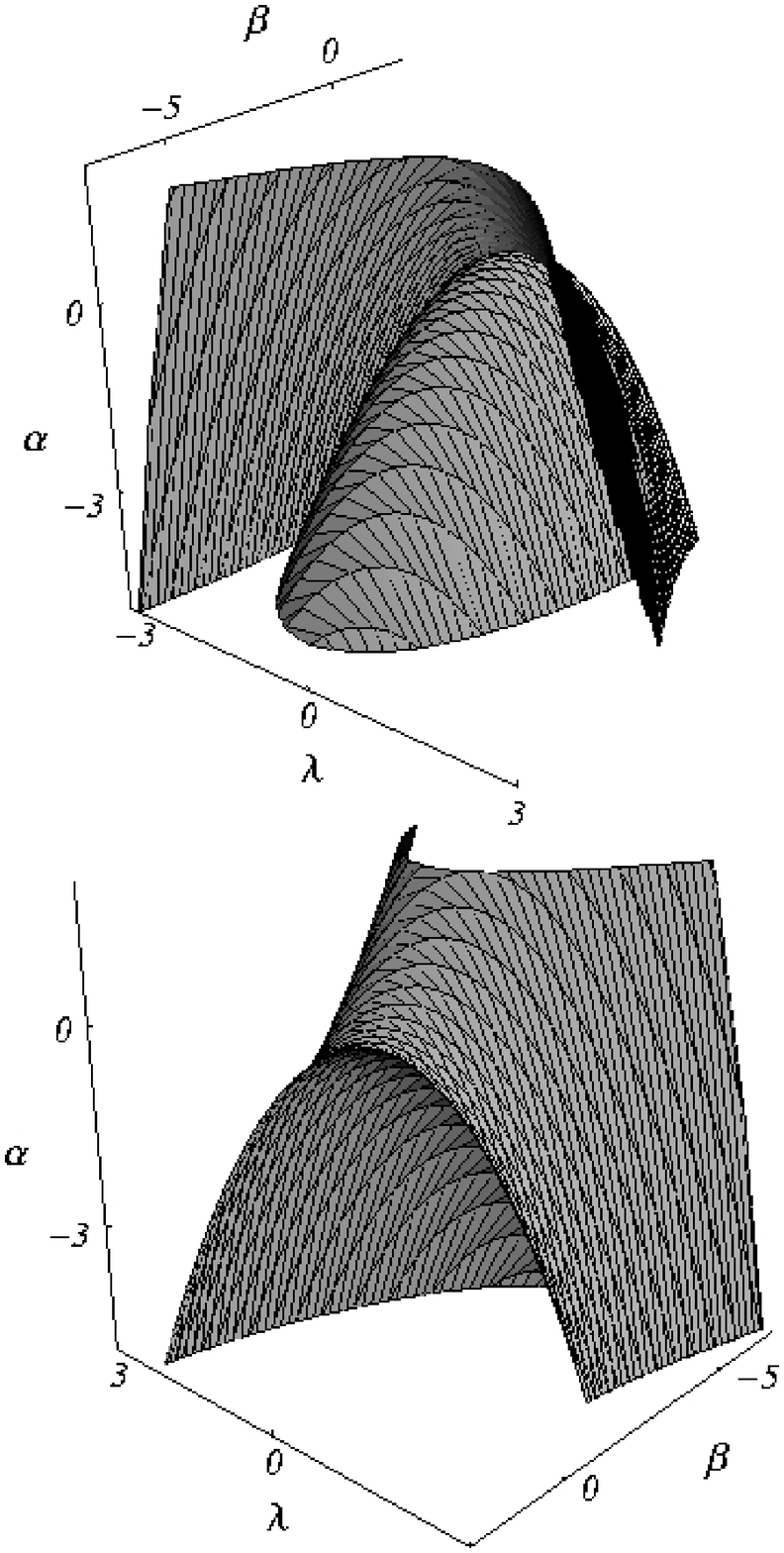,width=10cm}
\caption{\label{fig10}Front and back views of
$L_{\scriptscriptstyle TI}$. $P_{\scriptscriptstyle TI}$ occurs at (0,0,0).}
\end{figure}

One of the greatest strengths of singularity theory is that it defines criteria
by which two algebraic systems have equivalent solution sets.  
Singularity theory tells us that the CSTR problem is qualitatively equivalent
to the simplest universal unfolding of a bifurcation problem containing 
$P$, $H$, $T$, $I$, and $A$. This is designated $P_{\scriptscriptstyle TI}$:
\begin{equation}\tag{$P_{\scriptscriptstyle TI}$}
G(x,\lambda,\alpha,\beta)
=x^3+\lambda(\lambda-x)+ \alpha+\beta x.
\end{equation}
Evaluation of $T$ and $I$ embedded in $P_{\scriptscriptstyle TI}$ and $A$
(Appendix B) indicates that the 
bifurcation behaviour of the CSTR problem is completely 
encapsulated in the simpler problem 
$P_{\scriptscriptstyle TI}$. The limit point shell $L_{\scriptscriptstyle TI}$ 
of  $P_{\scriptscriptstyle TI}$
is viewed from two vantage points in figure~\ref{fig10}, and it 
appears to tell a different story. 

\begin{figure}\hspace*{1cm}
\psfig{file=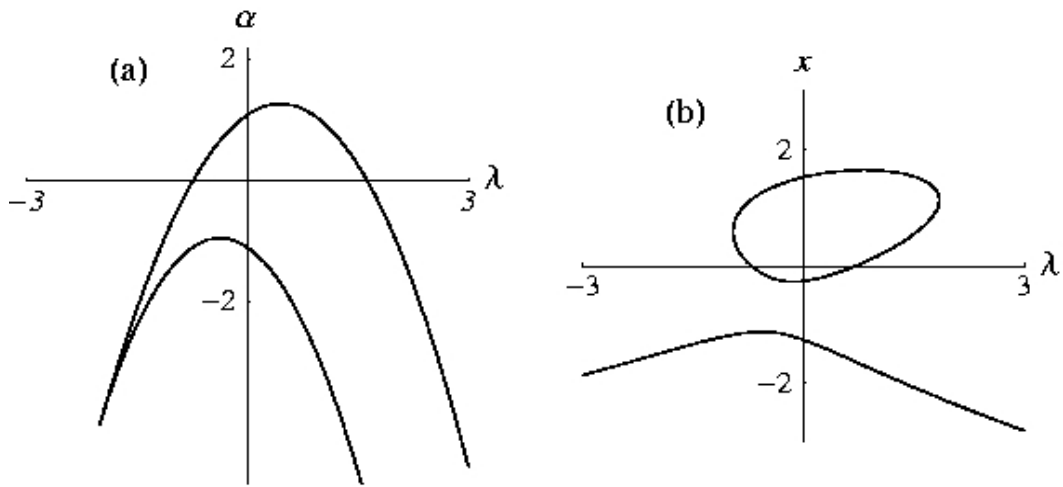,width=14cm}
\caption{\label{fig11}(a) A cross section of $L_{\scriptscriptstyle 
TI}$ at $\beta=-2$
gives the fold lines in the $\lambda,\alpha$ plane. (b) A cross section
at $\alpha=-0.5$ of the cross section at $\beta=-2$ gives one of the possible
bifurcation diagrams for $P_{\scriptscriptstyle TI}$.}
\end{figure}

All of the qualitatively distinct bifurcation diagrams of the CSTR can 
indeed be recovered from various slices of 
$L_{\scriptscriptstyle TI}$ --- for example, figure \ref{fig11} shows the
isolated branch of solutions that also exists in the CSTR. 
However, the qualitative 
differences between the two limit-point shells are quite striking. We could
not use $P_{\scriptscriptstyle TI}$ to predict the boundaries of 
multiplicity in the CSTR. 
\subsection{\label{lh} A tale of three pitchforks}
An important issue in the physics of magnetically confined plasmas 
is  the dramatic 
jump to an improved confinement r\'{e}gime, known as the L--H transition,
that occurs at critical values of
tunable parameters or internal system properties.  
A phenomenological model that described this critical behaviour 
\cite{Sugama:1995} was analysed by Ball and Dewar \cite{Ball:1999c} and
found to contain a partially unfolded pitchfork. The simplest universal 
unfolding for this bifurcation problem was derived as:  
\begin{align}
G(u,q,d,\alpha)=\frac{\left(dq-u^2\right)\left(b + au^{5/2}\right)}{u^{5/2}}
+\frac{q\left(u^2-dq\right)}{u^2}+\alpha,\tag{$P_{\scriptscriptstyle LH}$}
\end{align}
where $u$ is the state variable, $q$, $d$, and $\alpha$ are parameters and
$a$ and $b$ are numerical factors having the values 0.05 and 0.95 
respectively. In the model, $u$, $q$, and $d$ $> 0$. ({\it Note:} The dynamical
L-H model also contains Hopf bifurcations, which will be discussed elsewhere.)

Unusually, although there is a degenerate pitchfork in 
the original partially unfolded (or overdetermined)
problem with $\alpha=0$, the universal unfolding $P_{\scriptscriptstyle LH}$
 contains {\em two} pitchforks in the real, physical space. They 
are found
exactly by applying the conditions $P$ in table \ref{table1} 
to~$P_{\scriptscriptstyle LH}$:
\begin{align}
P_{1\scriptscriptstyle LH}: \quad&(u,q,d,\alpha)=
\left(\frac{b^{2/5}}{2^{4/5}a^{2/5}},\,
\frac{5b}{4\sqrt{u}},\,\frac{4u^{5/2}}{5b},\,0\right)\notag\\
P_{2\scriptscriptstyle LH}:  \quad&(u,q,d,\alpha)=\notag\\
&\left(\frac{\left(173+70\sqrt{6}\right)^{1/5}b^{2/5}}{22^{2/5}a^{2/5}},\,
\frac{5b\left(2+\sqrt{6}\right)}{16\sqrt{u}},\,
\frac{8u^{5/2}\left(2+\sqrt{6}\right)}{5b},\,
-\frac{5b\left(2+\sqrt{6}\right)}{2\sqrt{u}}\right)\notag
\end{align}
This introduces a formidable global aspect to what hitherto has been a purely
local focus on the structure of the limit-point shell around a {\em unique} 
pitchfork. A bifurcation analysis of $P_{\scriptscriptstyle LH}$
 and construction of
the limit-point shell is clearly
not for the faint-hearted, because it is difficult to imagine the
topological nightmares that must exist to connect 
$P_{1\scriptscriptstyle LH}$ and $P_{2\scriptscriptstyle LH}$. One can begin
by searching for lower-order and same-order bifurcations; the conditions for $T$ 
and $I$
yield:
\begin{align}
(u,q,d)&=u_{\scriptscriptstyle TI},\,
\left(\frac{5b\left(-5b+4\alpha\sqrt{u}\right)}{16\alpha u},\,
-\frac{16\alpha u^3}{25b^2}\right),\notag\\
\det{d^2G}&=-\frac{96\alpha^2}{25b^2u},\,\,
G_{uu}=\frac{25b^2/\alpha-20b\sqrt{u}-2\alpha u}{2u^3},\notag
\end{align} 
with $u_{\scriptscriptstyle TI}$ given by real roots of 
$16\alpha u^3a-b(-25b+24\alpha\sqrt{u})=0$. From the condition
$\det{d^2G}\neq 0$ we can infer that the bifurcation is transcritical, 
and since $u$, $q$, and $d$ $>0$
it is required that $\alpha <0$ at this point. The conditions for $H$
evaluate to expressions too complicated to be useful. 
There does not appear to be an asymmetric cusp in the problem. 

\begin{figure}\hspace*{2cm}
\psfig{file=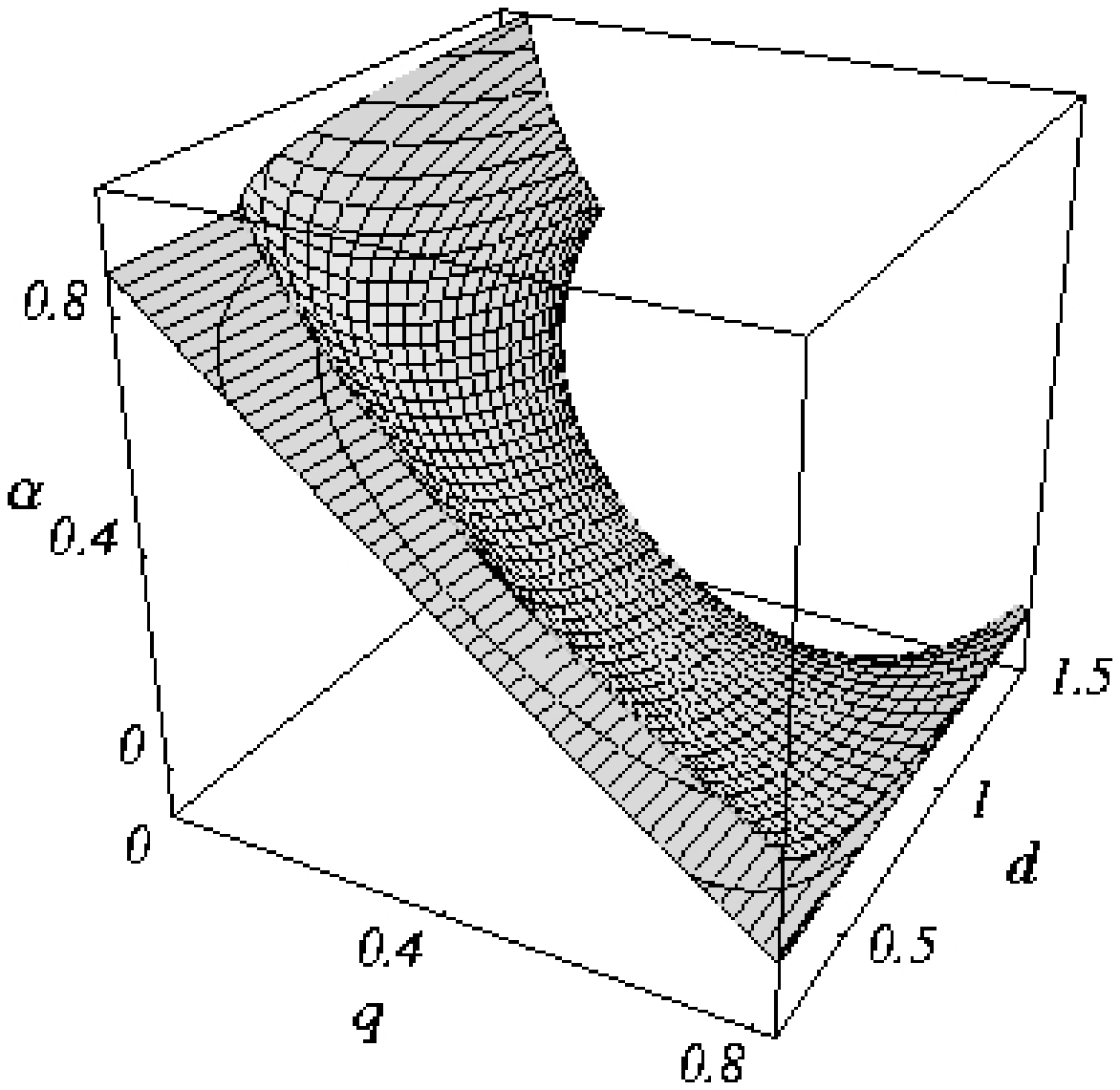,width=10cm}
\caption{\label{fig12}Part of $L_{\scriptscriptstyle LH}$ showing how 
$P_{1\scriptscriptstyle LH}$ and 
$P_{3\scriptscriptstyle LH}$ are connected.}
\end{figure}

\begin{figure}\hspace*{2cm}
\psfig{file=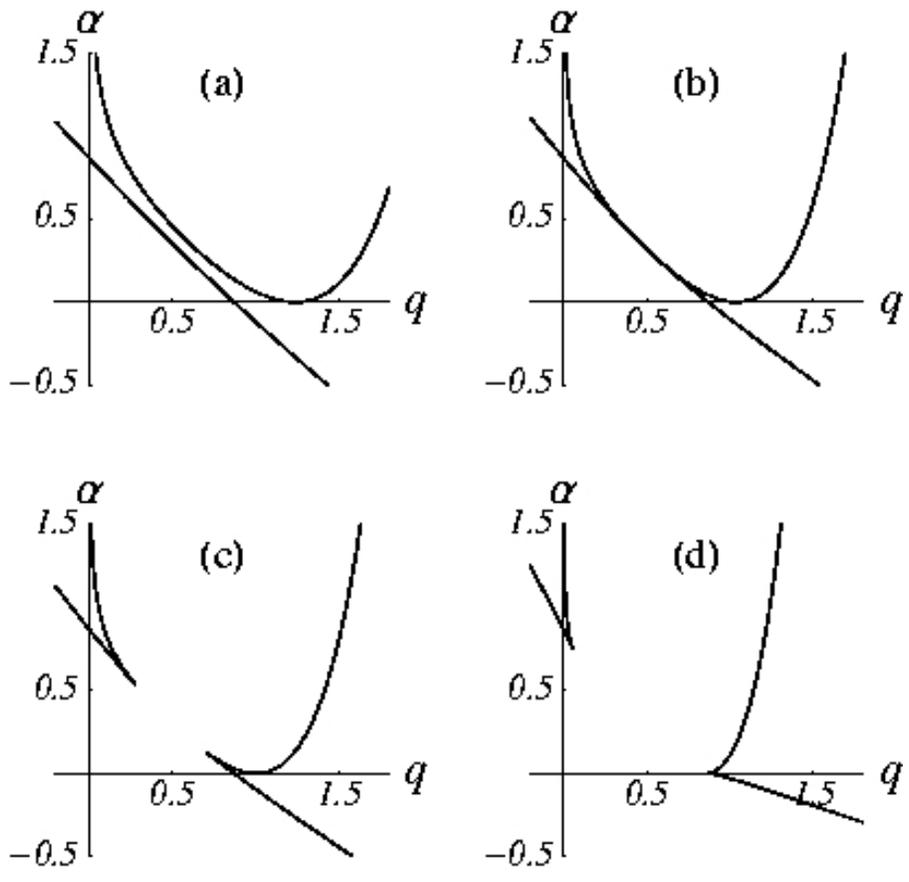,width=12cm}
\caption{\label{fig13}Successive slices of $L_{\scriptscriptstyle LH}$ 
illustrate the $E_2$ bifurcation. (a) $d=0.3$, (b) $d=0.8$, the 
critical value at $E_2$,
(c) $d=1.0$, (d) $d=4$, the critical value at $P_{1LH}$.}
\end{figure}

Another, related complication is the existence of a {\em third} pitchfork in 
the unphysical region $q<0$, $d<0$. It is given by
\begin{align}
P_{3\scriptscriptstyle LH}: \quad&(u,q,d,\alpha)=\notag\\
&\left(\frac{\left(173-70\sqrt{6}\right)^{1/5}b^{2/5}}{22^{2/5}a^{2/5}},\,
-\frac{5b\left(-2+\sqrt{6}\right)}{16\sqrt{u}},\,
-\frac{8u^{5/2}\left(-2+\sqrt{6}\right)}{5b},\,
\frac{5b\left(-2+\sqrt{6}\right)}{2\sqrt{u}}\right).\notag
\end{align}
$P_{3\scriptscriptstyle LH}$ is important because part of its limit-point 
manifold intrudes into the physical parameter region and is connected to the  
limit-point manifold of $P_{1\scriptscriptstyle LH}$ by a seam of hysteresis
points. This curved seam can be seen in figure \ref{fig12}, a
fragment of the limit-point shell of $P_{\scriptscriptstyle LH}$
 around the connection. In
figure \ref{fig13} a series of slices taken 
in the $q,\alpha$ plane illustrate how the connection occurs.

The putative condition for {\em \'{e}p\'{e}e \`{a} \'{e}p\'{e}e} contact 
of the two hysteresis points in figure~\ref{fig13},
a singularity designated $E_2$, appears to be
\begin{equation}\label{E2}
G=G_u=G_{uq}=G_{uu}=0,\,G_q>0,\,G_{uuu}<0.
\end{equation}
This condition yields a single exact $E_2$ point:
\begin{align}
E_2: \quad& \left(u,q,d,\alpha\right)=
\left(\frac{b^{2/5}}{22^{1/5}a^{2/5}},\,
\frac{5b}{8\sqrt{u}},\,
\frac{8u^{5/2}}{25b},\,
\frac{2b}{5\sqrt{u}}\right),\notag
\end{align}
at which $G_q=24/25$ and $G_{uuu}=-3b/u^{7/2}$. The conditions (\ref{E2}) for 
the $E_2$ singularity seem pathological, but they can be understood by
referring to the surface in figure \ref{fig12} and considering the
the singularities defined in table \ref{table1}. At the hysteresis points, $H$,
$G_q\neq 0$, and this condition holds at $E_2$ (since it is not a pitchfork).
However, another derivative must be zero at the union of two degenerate
points, this is the condition $G_{uq}=0$. 
\section{\label{s4}Discussion and conclusion}
The above analysis highlights the pitfalls of accepting qualitative 
equivalence of a bifurcation problem to a normal form as, in some sense,
a ``solution''. In the case of the CSTR problem and its normal form, 
the boundaries of multiplicity are profoundly different. The analysis of the
L--H problem also hints at the bizarre and interesting 
features that a limit-point shell can
have while still remaining continuous. Although 
$L_{\scriptscriptstyle LH}$ is locally equivalent to $L_p$ around each of the
 three pitchforks, 
the global definition of $L_{\scriptscriptstyle LH}$ involves at least one 
new singularity,
$E_2$.

This is largely an interpretive and exploratory work, investigating the 
local and global environment of the pitchfork through the limit-point shell
and using prototypic normal forms and real-world bifurcation problems as 
examples.
It turns out that 
bifurcation problems containing an organizing centre as simple as the 
pitchfork can have rather complex boundaries of multiplicity. For this reason,
and given the increasing availability of 3-dimensional computer 
visualization techniques, the limit-point shell has enormous potential as
a design and control aid for experimental dynamical systems. 
\subsubsection*{Acknowledgement:} I would like to thank Henry Gardner for 
helpful comments on this manuscript. 
\section*{Appendix}
\subsection*{\label{A}A.}
The conditions for $T$, $I$, and $A$ in table \ref{table1} are applied to
$P_c$. The condition for $P$ has been evaluated in 
\cite{Ball:1999a}. With $\ell$ fixed and nonzero, $T$ or $I$ points exist at 
\begin{align}
&\left(f,\theta,\varepsilon\right)=\notag\\
&\left(\frac{\ell^{1/3}u^{2/3}}{e^{-2/3u}-e^{-1/u}\ell^{1/3}u^{2/3}},\,
u-\frac{e^{-1/3u}u^{4/3}}{\ell^{1/3}},\,
\frac{e^{-1/3u}\ell^{1/3}\left(-1+e^{-1/3u}\ell^{1/3}u^{2/3}\right)^2}
{u^{2/3}}\right).\notag
\end{align}
(The implicit function theorem ensures that $u$ can in principle 
be given by these expressions.)
The nondegeneracy conditions at $T$ or $I$ points evaluate to
\begin{align}
G_{uu}&= u^{-8/3}e^{-1/3u}\ell^{2/3}\left(-1+2e^{1/3u}\ell^{1/3}u^{2/3}
-2u\right), \notag\\
\det d^2G&=u^{-8/3}e^{-2/3u}\ell^{2/3}\left(-1+e^{1/3u}
\ell^{1/3}u^{2/3}\right)^3\left(-1+3e^{1/3u}\ell^{1/3}u^{2/3}-4u\right).
\notag
\end{align}
These expressions may be positive, negative, or zero. Where $\det d^2G=0$ the
conditions $A$ in table \ref{table1} give an asymmetric cusp for $0<u<0.5$, at
\begin{align}
&\left(f,\theta,\varepsilon,\ell\right)=\notag\\
&\left(\frac{e^{-1/u}\left(-1-4u\right)}{-2+4u},\,
\frac{u\left(1+u\right)}{1+4u},\,
\frac{4}{27}\left(16u+\frac{1}{u^2}-12\right),\,
\frac{e^{-1/u}\left(1+4u\right)^3}{27u^2}\right).\notag
\end{align}
Since an asymmetric cusp is the point of coincidence of a transcritical point 
and an isola point, it may be inferred that at least one isola and one
transcritical bifurcation exists in the CSTR problem. 
\subsection*{B.}
The conditions in table \ref{table1} are applied to
$P_{\scriptscriptstyle TI}$. The pitchfork occurs at 
$(x,\lambda,\alpha,\beta)=(0,0,0,0)$, at which 
$G_{xxx}=6$, $G_{x\lambda}=-1$. $I$ occurs for 
\begin{align}
\left(x,\lambda,\alpha\right)&=
\left(\frac{1}{12}\left(1+B\right),\,
\frac{1}{24}\left(1+B\right),\,
\frac{1}{864}\left(1+B-
24\beta\left(3+2B\right)\right)\right),\notag\\
\mathrm{det }\,d^2G&=B,\quad
G_{xx}=\frac{1}{2}\left(1+B\right);\notag
\end{align}
$T$ occurs for 
\begin{align}
\left(x,\lambda,\alpha\right)&=
\left(\frac{1}{12}\left(1-B\right),\,
\frac{1}{24}\left(1-B\right),\,
\frac{1}{864}\left(1-B-
24\beta\left(3-2B\right)\right)\right),\notag\\
\mathrm{det }\,d^2G&=-B,\quad
G_{xx}=\frac{1}{2}\left(1-B\right);\notag
\end{align}
where$B=\sqrt{1-48\beta}.$ The asymmetric cusp occurs at 
$(x,\lambda,\alpha,\beta)=(1/12,1/24,-1/1728,1/48)$.
\renewcommand{\baselinestretch}{1}\normalsize

\newpage
\listoffigures
\listoftables
\newpage
\clearpage
\clearpage

 \begin{table}\hspace*{3cm} 
 \begin{tabular}{c|c c c c c}
 & $P$ & $H$ & $T$ & $I$ &$A$\\
 \hline
$G_\lambda$ & 0 & $\neq 0$ & 0 & 0 & 0\\
$G_{xx}$ & 0 & 0 & $\neq 0$ & $\neq 0$ &$\neq 0$\\
$G_{x\lambda}$&$\neq 0$&&&&\\
$G_{xxx}$&$\neq 0$&$\neq 0$&&&\\
det $d^2G$&&&$<0$&$>0$&0\\
$Q_3$&&&&&$\neq 0$\\
codimension&2&1&1&1&2\\ \hline
\end{tabular} 
\vspace*{1cm}\caption{\label{table1}Conditions on a bifurcation problem 
$G(x,\lambda,\alpha_i)$ for the pitchfork $P$, 
hysteresis $H$, transcritical $T$, isola $I$, and asymmetric cusp $A$
singularities.
The primary singularity conditions 
$G=G_x=0$ are assumed. The codimension indicated is the minimum number of
unfolding parameters (the $\alpha_i$) needed for a universal unfolding. $d^2G$
is the Hessian, or matrix of second partial derivatives. $Q_3$ is a third-order
directional derivative, explained in \cite{Golubitsky:1985}.}
\end{table}   

\begin{thebibliography}{10}

\bibitem{Kuznetsov:1998}
Y.~A. Kuznetsov.
\newblock {\em Elements of Applied Bifurcation Theory}.
\newblock Springer-Verlag, New York, 2nd edition, 1998.
\newblock Applied mathematical sciences vol. 112.

\bibitem{Broer:1999}
H.~W. Broer, I.~Hoveijn, M.~van Noort, and G.~Vegter.
\newblock The inverted pendulum: a singularity theory approach.
\newblock {\em Journal of Differential Equations}, 157:120--149, 1999.

\bibitem{Schmitt:1998}
J.~M. Schmitt and P.~V. Bayly.
\newblock Bifurcations in the mean angle of a horizontally shaken pendulum:
  analysis and experiment.
\newblock {\em Nonlinear Dynamics}, 15(1):1--14, 1998.

\bibitem{Johnson:1998}
R.~C. Johnson.
\newblock Unicycles and bifurcations.
\newblock {\em American Journal of Physics}, 66(7):589--592, 1998.

\bibitem{Yabuno:1998}
H.~Yabuno, Y.~Kurata, and N.~Aoshima.
\newblock Effect of coulomb damping on buckling of a two-rod system.
\newblock {\em Nonlinear Dynamics}, 15(3):207--224, 1998.

\bibitem{Golubitsky:1979}
M.~Golubitsky and D.~G. Schaeffer.
\newblock Imperfect bifurcation in the presence of symmetry.
\newblock {\em Communications in Mathematical Physics}, 67(3):205--232, 1979.

\bibitem{Jensen:1999}
C.~N. Jensen, M.~Golubitsky, and H.~True.
\newblock Symmetry, generic bifurcations, and mode interaction in nonlinear
  railway dynamics.
\newblock {\em International Journal of Bifurcation and Chaos},
  9(7):1321--1331, 1999.

\bibitem{Algaba:1999}
A.~Algaba, E.~Freire, E.~Gamero, and A.~J. Rodr{\'\i}guez-Luis.
\newblock {A three-parameter study of a degenerate case of the Hopf-pitchfork
  bifurcation}.
\newblock {\em Nonlinearity}, 12(4):1177--1206, 1999.

\bibitem{Gray:1990}
P.~Gray and S.~K. Scott.
\newblock {\em Chemical Oscillations and Instabilities}.
\newblock Oxford University Press, 1990.

\bibitem{Sugama:1995}
H.~Sugama and W.~Horton.
\newblock {L--H confinement mode dynamics in three-dimensional state space}.
\newblock {\em Plasma Phys. Control. Fusion}, 37:345--362, 1995.

\bibitem{Ball:1999c}
R.~Ball and R.~L. Dewar.
\newblock {Singularity theory study of overdetermination in models for L--H
  transitions}.
\newblock Preprint URL: http://xxx.lanl.gov/abs/math-ph/9908023 or
  http://xxx.adelaide.edu.au/abs/math-ph/9908023, 1999.

\bibitem{Juel:1997}
A.~Juel, A.~G. Darbyshire, and T.~Mullin.
\newblock {The effect of noise on pitchfork and Hopf bifurcations}.
\newblock {\em Proceedings of the Royal Society of London Series A},
  453(1967):2627--2647, 1997.

\bibitem{Stewart:1999}
D.~E. Stewart and R.~L. Dewar.
\newblock Dynamical systems: Physics and numerical analysis.
\newblock In {\em Complex Systems}. Cambridge University Press, Cambridge, U.K,
  1999.

\bibitem{Ball:1999a}
R.~Ball.
\newblock The origins and limits of thermal steady state multiplicity in the
  continuous stirred tank reactor.
\newblock {\em Proc. R. Soc. Lond. A}, 455:141--161, 1999.

\bibitem{Ball:1999b}
R.~Ball and B.~F. Gray.
\newblock {Thermal stabilization of chemical reactors. II. Bifurcation analysis
  of the Endex CSTR}.
\newblock {\em Proc. R. Soc. Lond. A}, 1999.
\newblock (to appear).

\bibitem{Whitney:1955}
H.~Whitney.
\newblock {On singularities of mappings of Euclidean spaces. I. Mappings of the
  plane into the plane}.
\newblock {\em Ann. Math. II}, 62:374--410, 1955.

\bibitem{Thom:1972}
R.~Thom.
\newblock {\em Stabilit{\'e} Structurelle et Morphog{\'e}n{\'e}se : Essai d'
  une Th{\'e}orie G{\'e}n{\'e}rale des Mod{\`e}les}.
\newblock W.A. Benjamin, Reading, Mass., 1972.

\bibitem{Arnold:1992}
V.~I. Arnol'd.
\newblock {\em Catastrophe Theory}.
\newblock Springer-Verlag, Berlin, 3rd edition, 1992.
\newblock Translated from the Russian by G. S. Wassermann, based on a
  translation by R.K. Thomas.

\bibitem{Golubitsky:1985}
M.~Golubitsky and D.~G. Schaeffer.
\newblock {\em Singularities and Groups in Bifurcation Theory}, volume~1.
\newblock Springer--Verlag, New York, 1985.

\bibitem{Golubitsky:1980}
M.~Golubitsky and B.~L. Keyfitz.
\newblock A qualitative study of the steady-state solutions for a continuous
  flow stirred tank reactor.
\newblock {\em SIAM J. Math. Anal.}, 11:316--339, 1980.

\end{thebibliography}
\end{document}